\documentclass[12pt]{amsart}
\usepackage{fullpage}
\usepackage{amsmath,amsfonts,amssymb}
\usepackage{tikz}
\usepackage{hyperref}
\usepackage{graphicx}
\hypersetup{colorlinks=true,citecolor=red,linkcolor=red,urlcolor=black}
\newtheorem{theorem}{Theorem}[section]
\newtheorem{lemma}[theorem]{Lemma}
\newtheorem{proposition}[theorem]{Proposition}
\newtheorem{corollary}[theorem]{Corollary}
\newtheorem{definition}[theorem]{Definition}
\newtheorem{remark}[theorem]{Remark}
\numberwithin{equation}{section}

\begin{document}
\author[Beslikas]{Athanasios Beslikas}
\address{Doctoral School of Exact and Natural Studies, Institute of Mathematics, Faculty of Mathematics and Computer Science, Jagiellonian University, Lojasiewicza 6, PL-30-348 Krak\'ow, Poland}
\email{athanasios.beslikas@doctoral.uj.edu.pl}
\author[Krajina]{Leonora Krajina}
\address{Doctoral School of Exact and Natural Studies, Institute of Mathematics, Faculty of Mathematics and Computer Science, Jagiellonian University, Lojasiewicza 6, PL-30-348 Krak\'ow, Poland}
\email{leonora.krajina@doctoral.uj.edu.pl}
\author[Torkinejad Ziarati]{Pouriya Torkinejad Ziarati}
\address{Institut de Math\'ematiques de Toulouse, Universit\'e Toulouse III--Paul Sabatier, Toulouse, France}
\email{pouriya.torkinejad\_ziarati@math.univ-toulouse.fr}
\title[Dirichlet-type membership of RIFs]{The Dirichlet space of the bidisc contains no singular rational inner functions.}
\begin{abstract}
We characterize the membership of rational inner functions in the Dirichlet-type spaces of the bidisc defined by Kaptano\u{g}lu, thereby resolving an open conjecture from past work of Bickel, Pascoe and Sola. In particular, we prove that no rational inner function with singularities belongs to the unweighted Dirichlet space $\mathfrak{D}(\mathbb{D}^2).$
\end{abstract}
\maketitle

\section{Introduction}
\subsection{Rational inner functions on the bidisc.}
Consider the bidisc $$\mathbb{D}^2:=\{z=(z_1,z_2)\in\mathbb{C}^2:|z_1|<1,\ |z_2|<1\}$$ and its \textit{distinguished boundary,} the two-torus, given by
$$\mathbb{T}^2:=\{\zeta=(\zeta_1,\zeta_2):|\zeta_1|=1,|\zeta_2|=1\}.$$
Let $p\in\mathbb{C}[z_1,z_2]$ be a polynomial of bidegree $\deg(p)=(m,n)\in\mathbb{N}_0^2.$ We say that such a polynomial is \textit{stable} if it does not have any zeros on the unit bidisc $\mathbb{D}^2.$ Furthermore, consider its reflection
$$\widetilde{p}(z_1,z_2):=z_1^mz_2^n\overline{p\left(\frac{1}{\overline{z_1}},\frac{1}{\overline{z_2}}\right)},\quad z=(z_1,z_2)\in \mathbb{D}^2.$$
By the classical theorem of Rudin, Stout and Pfister (see for example \cite{Rudin}), all rational inner functions (briefly RIFs) on the bidisc admit the representation
$$\varphi(z_1,z_2)=e^{i\theta} z_1^Mz_2^N\frac{\widetilde{p}(z_1,z_2)}{p(z_1,z_2)}$$
for $(z_1,z_2)\in\mathbb{D}^2$ and some $(M,N)\in\mathbb{N}_0^2$. Recall that a bounded holomorphic function $\varphi$ is called \textit{inner} if $|\varphi^*|=1$ almost everywhere on $\mathbb{T}^2,$ where $\varphi^*$ denotes its radial limit, given almost everywhere by $$\varphi^{*}(\zeta_1,\zeta_2)=\lim_{r\to1^{-}}\varphi(r\zeta_1,r\zeta_2),$$ for $(\zeta_1,\zeta_2)\in\mathbb{T}^2.$ 

In one complex variable, rational inner functions are precisely the finite Blaschke products, which extend holomorphically to a disc strictly larger than $\mathbb{D}.$ In higher dimensions, rational inner functions may have singularities on $\mathbb{T}^2.$ As an example, consider the RIF
$$\kappa(z_1,z_2)=\frac{2z_1z_2-z_1-z_2}{2-z_1-z_2},\quad (z_1,z_2)\in\mathbb{D}^2,$$ which has a singularity at $(1,1)\in\mathbb{T}^2.$

\subsection{The Kaptano\u{g}lu definition of Dirichlet-type spaces on the bidisc.}
Let $f(z_1,z_2)=\sum_{k,\ell\geq 0}a_{k\ell}z_1^kz_2^\ell$ be holomorphic in $\mathbb{D}^2,$ i.e. $f\in\mathcal{O}(\mathbb{D}^2,\mathbb{C}).$ For $\vec{\alpha}=(\alpha_1,\alpha_2)\in\mathbb{R}^2,$ the Dirichlet-type spaces $\mathfrak{D}_{\vec{\alpha}}(\mathbb{D}^2)$ are defined by the coefficient norm
$$\mathfrak{D}_{\vec{\alpha}}(\mathbb{D}^2)=\Bigg\{f\in\mathcal{O}(\mathbb{D}^2,\mathbb{C}):\sum_{k,\ell\in\mathbb{Z_+}}(k+1)^{\alpha_1}(\ell+1)^{\alpha_2}|a_{k\ell}|^2<+\infty\Bigg\}.$$
For $\alpha_j\in(0,2),$ one has the equivalent integral characterization
$$\mathfrak{D}_{\vec{\alpha}}(\mathbb{D}^2):=\Biggl\{ f\in \mathcal{O}(\mathbb{D}^2,\mathbb{C}):\int_{\mathbb{D}^2}\left|\frac{\partial^2(z_1z_2f(z_1,z_2))}{\partial z_1 \partial z_2}\right|^2dA_{\alpha_1}(z_1)dA_{\alpha_2}(z_2)<\infty\Biggr\},$$ 
for a pair $\vec{\alpha}=(\alpha_1,\alpha_2)\in(0,2)^2,$ where
$dA_{\alpha_j}(z_j)=(1-|z_j|^2)^{1-\alpha_j}dA(z_j)$ and $dA(z_j)=dx_jdy_j/\pi$ are normalized Lebesgue measures on the unit disc.
We also write $\mathfrak{D}_{\alpha}=\mathfrak{D}_{(\alpha,\alpha)}$ and $\mathfrak{D}=\mathfrak{D}_{(1,1)}.$ These Dirichlet-type spaces on the bidisc were introduced by Kaptano\u{g}lu in \cite{PJM}.

The membership problem for RIFs was investigated in \cite{PLMS}, and the following conjecture was stated explicitly in \cite{Studia}.

\begin{itemize}
\item \textbf{Conjecture 1. }\textit{Let $\varphi=\frac{\tilde{p}}{p}$ be a rational inner function on the bidisc with at least one singularity. Then $\varphi\notin\mathfrak{D}(\mathbb{D}^2).$}
\end{itemize}
We resolve this conjecture and obtain the following stronger characterization.
\begin{theorem}\label{MainTheorem}
Let $\varphi$ be a RIF on $\mathbb{D}^2$ with at least one singularity on $\mathbb{T}^2,$ and let $K$ denote its contact order. Then, for every $\alpha\in\mathbb{R},$
$$\varphi\in\mathfrak{D}_{\alpha}(\mathbb{D}^2)\quad\iff\quad\alpha<\frac12+\frac{1}{2K}.$$
\end{theorem}
Here $K\ge2$ denotes the \textit{contact order} of $\varphi,$ an algebro-geometric quantity that we shall explain shortly in the next section.

As a straighforward consequence, we obtain the affirmative answer to this conjecture, that is
\begin{corollary}
Let $\varphi$ be a RIF on $\mathbb{D}^2$ with at least one singularity on $\mathbb{T}^2.$ Then $\varphi\notin\mathfrak{D}(\mathbb{D}^2).$
\end{corollary}
\section{Background}
\subsection{Contact order}  Consider a RIF
$\varphi=\widetilde{p}/p$ and assume that this function
has a singularity at $\tau=(\tau_1,\tau_2)\in\mathbb{T}^2.$
In other words, we assume that $p(\tau_1,\tau_2)=0$
and consequently that implies that $\widetilde{p}(\tau_1,\tau_2)=0$ too. 

For each $\zeta_2\in\mathbb{T},$ the function
$\mathbb{D}\ni z_1\mapsto\varphi_{\zeta_2}(z_1)
=\varphi(z_1,\zeta_2)$ is a finite Blaschke product,
after cancelling the possible common factors of the
numerator and the denominator. Denote its zero set on $\mathbb D$ by $\mathcal{Z}_{\zeta_2}=\{a_1(\zeta_2),\ldots,
a_N(\zeta_2)\}\subset\mathbb{D}$.
As is explained in \cite{PLMS}, the distance of the
zero set of $\varphi_{\zeta_2}$ to the unit circle
$\mathbb{T}$ is given by the quantity
$$\epsilon(\varphi,\zeta_2)=
\min_{1\leq j\leq N}\{1-|a_j(\zeta_2)|\}.$$
As $\zeta_2\to\tau_2,$ some of the zeros of the finite
Blaschke products $\varphi_{\zeta_2}$ tend to $\tau_1.$
A natural question is at what ``speed" this phenomenon
occurs at. The answer to this question is given in the
following theorem.

\begin{theorem}(Theorem 3.5 of \cite{PLMS})
\label{def:loccont}
Let $\varphi=\widetilde{p}/p$ be an RIF with a singularity
on $\mathbb{T}^2$ with $z_2$-coordinate $\tau_2.$
Then there exists a rational number $K>0$ such that
$$\epsilon(\varphi,\zeta_2)\asymp|\tau_2-\zeta_2|^K,$$
for all $\zeta_2\in\mathbb{T}$ in a sufficiently small
punctured neighborhood of $\tau_2.$ The number $K$
is called the $(z_1,\tau_2)$-contact order of $\varphi.$
\end{theorem}

This result allows one to formally define the contact
order of an RIF, even in the case of multiple singularities.
More concretely,

\begin{definition}(Definition 3.6 of \cite{PLMS})
Let $\varphi=\widetilde{p}/p$ be a RIF and let
$\tau_{2,1},\ldots,\tau_{2,J}$ denote the distinct
$z_2$-coordinates of the singularities of $\varphi.$
Then the $z_1$-contact order of $\varphi$ is the number
$K^1$ defined by
$$K^1:=\max\{(z_1,\tau_{2,j})\text{-contact order of }
\varphi:1\leq j\leq J\}.$$
\end{definition}

The definition of the $z_2$-contact order $K^2$ of
$\varphi$ is set up in an analogous manner. It is known that $K^1=K^2,$
and that the common contact order $K$ is an even integer; see
\cite{Pisa,Polonici} for further details.
\subsection{Previous results}
In \cite[Corollary 2.2]{Studia}, the first named author and A. Sola proved the following result.
\begin{theorem} Let $\varphi$ be a RIF with at least one singularity on $\mathbb{T}^2,$ and let $K$ denote its contact order. Then, for $0<\alpha<2,$ we have $\varphi\in\mathfrak{D}_{(\alpha,0)}(\mathbb{D}^2)\cap\mathfrak{D}_{(0,\alpha)}(\mathbb{D}^2)$ if and only if $\alpha<1+\frac{1}{K}.$
\label{anisotropicchar}
\end{theorem}
As one observes, a complete characterization was obtained for the intersection $\mathfrak{D}_{(\alpha,0)}\cap\mathfrak{D}_{(0,\alpha)}$ in terms of contact order. However, the corresponding necessity condition for $\mathfrak{D}_{\alpha}$ was left open. In the present note, we observe that for an RIF, the membership in $\mathfrak{D}_{(\alpha_1,\alpha_2)}$ is actually equivalent to the membership in $\mathfrak{D}_{(\alpha_1+\alpha_2,0)}(\mathbb{D}^2),$ and equivalent to membership in $\mathfrak{D}_{(0,\alpha_1+\alpha_2)}(\mathbb{D}^2).$ 
Taking $\alpha_1=\alpha_2$ provides the missing necessity condition in Theorem \ref{MainTheorem}.

\section{Holomorphic extension and coefficient estimates}

The central point is a simple observation concerning stable polynomials. We write $\mathcal{Z}_p$ for the zero set of $p.$

\begin{lemma}\label{ExtensionLemma}
Let $\varphi=\widetilde{p}/p$ be a reduced RIF with a singularity at $\tau\in\mathbb{T}^2.$ Then there exist $R>1$ and $0<r<1$ such that $\varphi$ extends holomorphically to neighborhoods of both polydiscs $R\mathbb{D}\times r\mathbb{D}$ and $r\mathbb{D}\times R\mathbb{D}.$
\end{lemma}

\begin{proof}
First, by stability and our assumptions, $\mathcal{Z}_p\cap(\mathbb{T}\times\mathbb{D})=\mathcal{Z}_p\cap(\mathbb{D}\times\mathbb{T})=\varnothing.$
Fix now $0<r<1.$ Again by stability, and the preceding observation, $p$ is nonvanishing on the compact sets $\overline{\mathbb{D}}\times r\overline{\mathbb{D}}$ and $r\overline{\mathbb{D}}\times\overline{\mathbb{D}}.$ By Hartogs' theorem, we can extend the function to the Hartogs figure, that is, there exists $R>1$ such that $p$ is nonvanishing on $R\overline{\mathbb{D}}\times r\overline{\mathbb{D}}$ and $r\overline{\mathbb{D}}\times R\overline{\mathbb{D}}.$ The assertion follows.
\end{proof}

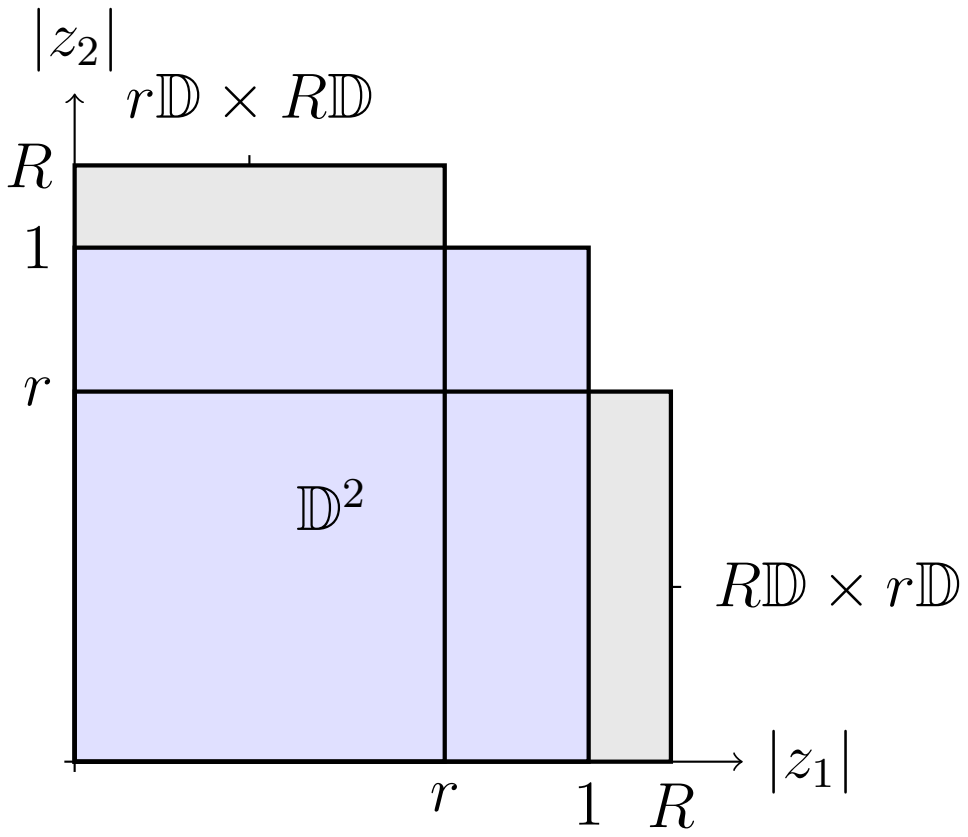
\begin{figure}[ht]
\begin{center}
\begin{tikzpicture}[scale=3.4]
\def\r{0.72}
\def\R{1.16}

\fill[gray!18] (0,0) rectangle (\R,\r);
\fill[gray!18] (0,0) rectangle (\r,\R);
\fill[blue!12] (0,0) rectangle (1,1);

\draw[thick] (0,0) rectangle (1,1);
\draw[thick] (0,0) rectangle (\R,\r);
\draw[thick] (0,0) rectangle (\r,\R);
\draw[dashed] (1,0) -- (1,1);
\draw[dashed] (0,1) -- (1,1);

\draw[->] (-0.02,0) -- (1.30,0) node[right] {$|z_1|$};
\draw[->] (0,-0.02) -- (0,1.30) node[above] {$|z_2|$};

\node at (0.50,0.50) {$\mathbb{D}^2$};

\node[anchor=west] at (1.20,0.34)
    {$R\mathbb{D}\times r\mathbb{D}$};
\draw (\R,0.34) -- (1.18,0.34);

\node[anchor=south] at (0.34,1.20)
    {$r\mathbb{D}\times R\mathbb{D}$};
\draw (0.34,\R) -- (0.34,1.18);

\node[below] at (1,0) {$1$};
\node[left] at (0,1) {$1$};
\node[below] at (\R,0) {$R$};
\node[left] at (0,\R) {$R$};
\node[below] at (\r,0) {$r$};
\node[left] at (0,\r) {$r$};

\end{tikzpicture}
\end{center}
\caption{Schematic representation in the $(|z_1|,|z_2|)$-plane
of the bidisc $\mathbb{D}^2$ (shaded blue) and the polydiscs
$R\mathbb{D}\times r\mathbb{D}$ and
$r\mathbb{D}\times R\mathbb{D}$, where $0<r<1<R$ (shaded gray).}
\end{figure}

The previous lemma provides an estimate on the Taylor coefficients of $\varphi$ which will allow us to compare the various Dirichlet-type norms.

\begin{proposition}\label{CoefficientProposition}
Let $\varphi(z_1,z_2)=\sum_{k,\ell\geq 0}a_{k\ell}z_1^kz_2^\ell$ be holomorphic in neighborhoods of $R\overline{\mathbb{D}}\times r\overline{\mathbb{D}}$ and $r\overline{\mathbb{D}}\times R\overline{\mathbb{D}},$ for some $R>1$ and $0<r<1.$ Then, for every $\alpha_1,\alpha_2\in\mathbb{R},$
$$\varphi\in\mathfrak{D}_{(\alpha_1,\alpha_2)}(\mathbb{D}^2)\iff \varphi\in\mathfrak{D}_{(\alpha_1+\alpha_2,0)}(\mathbb{D}^2)\iff \varphi\in\mathfrak{D}_{(0,\alpha_1+\alpha_2)}(\mathbb{D}^2).$$
\end{proposition}

\begin{proof}
By Cauchy's inequalities, there exists a constant $C>0$ such that
\begin{equation}\label{CauchyEstimates}
|a_{k\ell}|\leq\frac{C}{R^kr^\ell}\qquad \text{ and } \qquad |a_{k\ell}|\leq\frac{C}{r^kR^\ell},\quad k,\ell\geq0.
\end{equation}
Choose now $c>1$ sufficiently large so that $\rho=Rr^{1/c}>1.$ We divide the indices into three regions (see Figure \ref{regions}): 
\begin{itemize}
\item$S_1=\{(k,\ell):k+1>c(\ell+1)\},$
\item$S_2=\{(k,\ell):\ell+1>c(k+1)\}$ and 
\item $S_3=\mathbb{N}_0^2\setminus(S_1\cup S_2).$
\end{itemize}

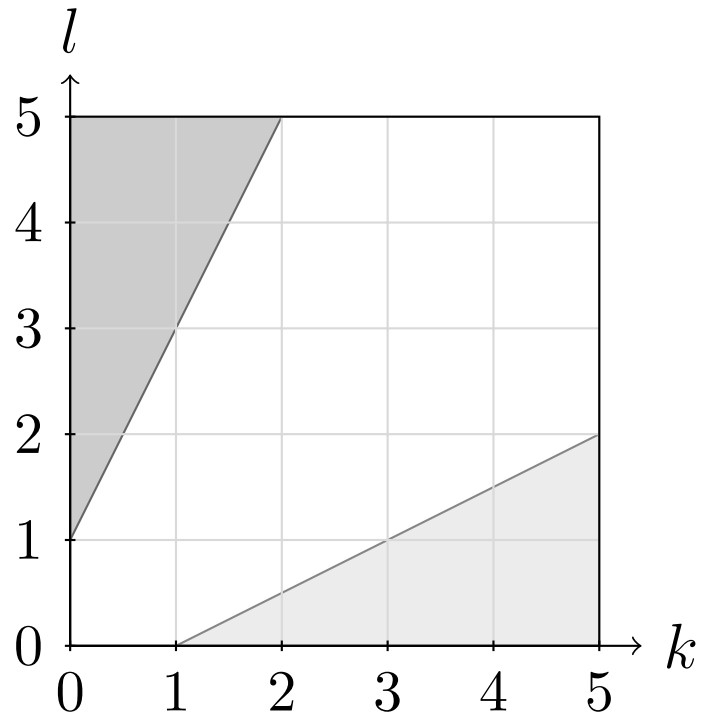
\begin{figure}[ht]
\begin{center}
\begin{tikzpicture}[scale=0.7]
  \fill[lightgray!30]   (1,0) -- (5,0) -- (5,2) -- cycle;
  \fill[gray!40] (0,1) -- (0,5) -- (2,5) -- cycle;
  \draw[lightgray!70!black]  (1,0) -- (5,2);
  \draw[gray!80!black] (0,1) -- (2,5);
  \draw[step=1, gray!30] (0,0) grid (5,5);
  \draw[->] (0,0) -- (5.4,0) node[right] {$k$};
  \draw[->] (0,0) -- (0,5.4) node[above] {$l$};
  \draw (0,0) rectangle (5,5);
  \foreach\x in {0,...,5} {
    \draw (\x,0.05) -- (\x,-0.05) node[below] {\small \x};
    \draw (0.05,\x) -- (-0.05,\x) node[left]  {\small \x};
}
\end{tikzpicture}
\end{center}
\caption{Example of how (filled out in $\mathbb{R}_+)$ regions $S_1$ (in light gray), $S_2$ (in gray) and $S_3$ (in white) might look like, for $c=2$, with $k$ and $\ell$ on the $x$- and $y$-axes, respectively.}
\label{regions}
\end{figure}

Observe that for $(k,\ell)\in S_1$ we have $\ell<k/c.$ Consequently, the first estimate in \eqref{CauchyEstimates} gives $|a_{k\ell}|\leq C\rho^{-k}.$ Since $\ell\leq k$ on $S_1,$ for every $\alpha_1,\alpha_2\in\mathbb{R}$ one obtains
$$\sum_{(k,\ell)\in S_1}(k+1)^{\alpha_1}(\ell+1)^{\alpha_2}|a_{k\ell}|^2\lesssim\sum_{k=0}^{\infty}(k+1)^{|\alpha_1|+|\alpha_2|+1}\rho^{-2k}<\infty.$$
The same argument, using the second estimate in \eqref{CauchyEstimates}, shows that the corresponding weighted sums over $S_2$ are finite.

Finally, for $(k,\ell)\in S_3$ we have
$$\frac1c\leq\frac{k+1}{\ell+1}\leq c.$$
Summarizing, we receive
$$\sum_{(k,\ell)\in S_3}(k+1)^{\alpha_1}(\ell+1)^{\alpha_2}|a_{k\ell}|^2\asymp\sum_{k=0}^{\infty}(k+1)^{\alpha_1+\alpha_2}|a_{k\ell}|^2= \|\varphi\|^2_{\mathfrak{D}_{(\alpha_1+\alpha_2,0)}(\mathbb{D}^2)}$$
and
$$\sum_{(k,\ell)\in S_3}(k+1)^{\alpha_1}(\ell+1)^{\alpha_2}|a_{k\ell}|^2\asymp\sum_{k=0}^{\infty}(\ell+1)^{\alpha_1+\alpha_2}|a_{k\ell}|^2= \|\varphi\|^2_{\mathfrak{D}_{(0,\alpha_1+\alpha_2)}(\mathbb{D}^2)}.$$
Therefore, for an RIF $\varphi$, the norm $\|\varphi\|^2_{\mathfrak{D}_{(\alpha_1,\alpha_2)}(\mathbb{D}^2)}$ is finite if and only if $\|\varphi\|_{\mathfrak{D}_{(\alpha_1+\alpha_2,0)}(\mathbb{D}^2)}$ is finite. This gives the desired equivalence of Dirichlet-type memberships for $\varphi$ on the bidisk and completes the proof.
\end{proof}

\begin{remark}
The conclusion of Proposition \ref{CoefficientProposition} does not require $\varphi$ to be a RIF, yet for our purposes it can be applied for such functions.
\end{remark}

\section{Proof of the main result}

\begin{proof}[Proof of Theorem \ref{MainTheorem}]
By Lemma \ref{ExtensionLemma} and Proposition \ref{CoefficientProposition}, for every RIF $\varphi$ and every $s\in\mathbb{R},$ we have
\begin{equation}\label{MembershipEquivalence}
\varphi\in\mathfrak{D}_{(s,0)}(\mathbb{D}^2)\iff\varphi\in\mathfrak{D}_{(0,s)}(\mathbb{D}^2)\iff\varphi\in\mathfrak{D}_{(s,0)}(\mathbb{D}^2)\cap\mathfrak{D}_{(0,s)}(\mathbb{D}^2).
\end{equation}
After setting $s=2\alpha,$ Proposition \ref{CoefficientProposition} and \cite[Corollary 2.2]{Studia} give us the following for $0<\alpha<1$
$$\varphi\in\mathfrak{D}_{\alpha}(\mathbb{D}^2)\iff
\varphi\in\mathfrak{D}_{(2\alpha,0)}(\mathbb{D}^2)\cap\mathfrak{D}_{(0,2\alpha)}(\mathbb{D}^2)
\iff 2\alpha<1+\frac{1}{K},$$
which is equivalent to $\alpha<\frac12+\frac{1}{2K}.$
\end{proof}

\begin{corollary}
Let $\varphi$ be a RIF on $\mathbb{D}^2$ with at least one singularity on $\mathbb{T}^2.$ Then $\varphi\notin\mathfrak{D}(\mathbb{D}^2).$
\end{corollary}

\begin{proof}
Since $K\geq2,$ we have $2\geq1+1/K.$ The result follows after applying Theorem \ref{MainTheorem} with $\alpha=1.$
\end{proof}

\section{Acknowledgments}
We express our gratitude to Alan Sola for sharing the problem with us and encouraging us to write up the solution in this paper. Moreover, we also thank Prof. Pascal Thomas for the useful mathematical discussions. 
\section{Financial Support}
The second named author is financially supported by the National Science Centre Programme K/NCN/000315 under the research project OPUS 26.
    
\section{AI disclosure} The mathematical arguments presented in this paper are purely made out of humans, as well as the write up. We used Chat GPT 6.1. Sol for confirming the correctness of our proof as it seemed too good to be true, as well as for the correction of grammatical mistakes and Figure 1.

\end{document}